\documentclass[12pt]{article}
 \usepackage{amsmath,amsfonts,theorem}
 \numberwithin{equation}{section}

 \newtheorem{prop}{Proposition}[section]

 \newtheorem{cor}{Corollary}[section]

 \theorembodyfont{\rmfamily}

 \newcommand{\qed}{\ifhmode\unskip\nobreak\fi\quad\ensuremath\square}

 \newcommand{\sS}{\mathcal S}
 
 \newcommand{\al}{\alpha}
 \newcommand{\be}{\beta}
 \newcommand{\de}{\delta}

 \newcommand{\ga}{\gamma}

 \newcommand{\De}{\Delta}

 \newcommand{\la}{\lambda}

 \newcommand{\PP}{\mathbb P}
 \newcommand{\C}{\mathbb C}

 \newcommand{\R}{\mathbb R}

 \newcommand{\diag}{\operatorname{diag}}

\begin{document}

 \title{Twist tori and pseudo toric structures}
 \markright{\hfill Twist tori... \quad}

 \author{Nikolay A. Tyurin\footnote{BLTPh JINR (Dubna) and MGU PS (MIIT).  This work
 was partially supported by RFBR, grants
N 08-01-00095, N 08-01-00392a}}
\date{}
 \maketitle

\begin{center}
 {\em                                  }
 \end{center}
 \bigskip

\section*{Introduction}

Twist tori are examples of exotic monotone lagrangian tori, presented in [1].
This tree of examples grew up over the first one --- the torus $\Theta \in
\R^4$, constructured in [2] and [3]. On the other hand, in [4] and [5] we
proposed a new structure  which generalizes the notion of toric structure. One
calls this generalization pseudo toric structure, and several examples were
given which show that certain toric symplectic manifolds can carry the structre
and that certain non toric symplectic manifolds do the same.

Below we show that any twist torus $\Theta^k \subset \R^{2k+2}$, defined in
[1], can be constructed via pseudo toric considerations. Due to this one can
explicitly show that every $\Theta^k \subset \R^{2k+2}$ is displaceable.

\section{Pseudotoric structure on $\C^{k+1}$ and twist torus $\Theta^k$}

Consider $\R^{2k+2}$ with  standard complex and symplectic structures and
denote the endowed space as $\C^{k+1}$. Choose a complex coordinate system
$(z_1, ..., z_{k+1}$ on it. Then we have a fibration
$$
\psi: \C^{k+1} \to \C
$$
given by the formula
$$
(z_1, ..., z_{k+1}) \mapsto a = z_1 \cdot ... \cdot z_{k+1} \in \C,
$$
and we denote the target space of the map as $\C_a$. The fibers of the map
$\psi$ are hypersurfaces in $\C^{k+1}$ defined by the equation
$$
z_1 \cdot ... \cdot z_{k+1} = a;
$$
they are smooth except for $a = 0$, and the last one is formed by $k+1$ basis
hyperplanes.

The standard complex torus action of $T^{k+1}$ on $\C^{k+1}$ can be reduced to
a subtorus $T^k_0$ whose action preserves the fibers of $\psi$. We can describe
the corresponding Hamiltonian action of a real subtorus $T^k_{\R}$ explicitly:
take a set of self adjoint operators $A_1, ..., A_k$ given in the fixed
coordinate system by diagonal matrices such that
$$
A_i = \diag(\la_1, ..., \la_{k+1}), \vert  \la_i = 1, \la_{i+1} = -1,
$$
and $\la_j =0$ for the resting $\la$'s; consider the functions
$$
F_i = <A_i \psi, \psi>_h,
$$
which are known as {\it expectation values} in Quantum Mechanics. Then it's not
hard to see that:

--- $\{F_1, ..., F_k\}$ is a complete set of commuting functions;

--- for each $F_i$ its Hamiltonian action preserves each fiber of $\psi$.

The structure of this kind was called  {\it a pseudo toric structure} (the
detailed definition can be found in [5]).

It's well known that each smooth fiber of $N_a = \psi^{-1}(a)$ endowed with the
function set $\{ F_1^a, ..., F^a_k \}$, where
$$
F_i^a = F_i|_{N_a},
$$
 is a completely integrable
system, see [6]. If we fix a couple $(c_1, ..., c_k)$ of non critical values
for the functions $F^a_1, ..., F^a_k$ on a smooth fiber $N_a$ (so $ -1 < c_i <
1$) then the mutual level set
$$
S^a_{(c_1, ..., c_k)} = \{ F^a_i = c_i \} \subset N_a
$$
is a smooth lagrangian torus for any $a \neq 0$.

{\bf Remark.} Note that since $F_1, ..., F_k$ are global functions on
$\C^{k+1}$ then we can fix the values $(c_1, ..., c_k)$ universally for all
fibers $N_a$ and all functions $(F^a_1, ..., F^a_k)$. On the other hand the non
critical values for $F_i$ and $F^a_i, a \neq 0,$ are different: for the last
function set certain $c_i$s can equal to zero. This happens since the critical
points for $F_1, ..., F_k$ with the zero critical values lie on the singular
fiber $N_0$.

Now we have the following

\begin{prop} Fixing a couple $(c_1, ..., c_k), - 1 < c_i, < 1,$ of non critical
 values for $F^a_1, ..., F^a_k$
and a smooth loop $\ga \subset \C_a$ which doesn't pass through the origin, one
defines a lagrangian torus in $\C^{k+1}$.
\end{prop}

{\it Proof.} For all $a \in \ga \subset \C_a$ consider the corresponding tori
$S^a_{(c_1, ..., c_k)}$ in the fibers of $\psi$ and combine all the tori
$$
S_{\ga, (c_1, ..., c_k)} = \bigcup_{a \in \ga} S^a_{(c_1, ..., c_n)};
$$
obviously it gives us a smooth torus $S_{\ga, (c_1, ..., c_n)} \in \C^{k+1}$.

Let's show that this torus is lagrangian. The tangent bundle $TS_{\ga, (c_1,
..., c_2)}$ is spanned by $k+1$ non vanishing vector fields:

--- the Hamiltonian vector fields $X_{F_i}$ of the functions $F_1, ..., F_k$,

--- and one additional vector field $\nabla_{\psi} \rho_{\ga}$ which is given by the
 lift, generated by the symplectic connection $\nabla_{\psi}$ (note that the fibers of $\psi$ are
symplectic), of a non vanishing tangent vector field $\rho_{\ga}$ on the loop
$\ga$.

Since the vector field $\nabla_{\psi} \rho_{\ga}$ is symplectically orthogonal
to the fibers of $\psi$, the tangent space $T_p S_{\ga, (c_1, ..., c_n)}$ is
lagrangian. This ends the proof.

Following D. Auroux, see [7] (and [4]), we distinguish two cases:

1) (the standard type)  the case when $\ga$ is not contractible in $\C_a^* =
\C_a - \{ 0 \}$;

2) (the Chekanov type) the case when $\ga$ is contractible in $\C^*_a$.

The reason is the following. Let's take $(c_1 = ... = c_k = 0)$, and $\ga \in
\C_a$ defined by the condition $\ga_0 = e^{i \phi}$. Then it's not hard to see
that $S_{\ga_0, (0, ..., 0)}$ is a standard torus in $\C^{k+1}$.

On the other hand, the parameter space $\C_a$ covers the diagonal complex line
$\Delta \subset \C^{k+1}$ given by the condition
$$
z_1 = ... = z_{k+1},
$$
so
$$
z^{n+1}: \De \to \C_a.
$$
 This means that in the loop space for $\C_a$ one has a component which
consists of loops of the following type: take the sector
$$
\sS(k) = \{r e^{i \phi}|0 < \phi < \frac{2\pi}{k+1} \} \in \De
$$
and consider the loops in $\C_a$ which are given by loops in $\De$, decribed by
conditions at the end of page 1 in [1], and then mapped to $\C_a$ by the taking
of $n+1$th powers. Let's fix a smooth loop $\ga \subset \De$ lies in the sector
$\sS (k) \cap D^2(k+1 + \epsilon)$ as required in [1]. For each point of $\ga$
take the $k+1$th power of the corresponding complex number. Then we get a
smooth loop $\ga_1 \in \C_a$.  We have the following

\begin{prop} In the situation of the Proposition 1.1 let's take $c_1 = ... c_k =0$
and loop $\ga_1$  described above. Then the corresponding lagrangian torus
$S_{\ga_1, (0,..., 0)}$ concides with the twist torus $\Theta^k$ from [1].
\end{prop}

{\it Proof} is given just by direct calculation.

Thus the constructions coming with the pseudo toric structures cover and
generalize the construction of twist tori, presented in [1]. Indeed, to get the
twist torus $\Theta^k$ we must take the loops in $\C_a$ of rather special type
while the construction works for any smooth loop.

Pseudo toric considerations hint a way how certain Hamiltonian isotopies of
$\C^{k+1}$ can be explicitly constructed. We discuss it in the next section.

\section{Displaceability}

Consider again the same situation as in Proposition 1 above so when a couple of
values $(c_1, ..., c_k)$ and a loop $\ga \in \C^a$ defines a lagrangian torus
in $\C^{k+1}$. Suppose that $\ga$ is contractible (so the torus is of the
Chekanov type). Then one can find a smooth function $f$ on $\C_a$ which is

(1) identically zero on a small disc $D^2(\de) = \{ \vert w \vert < \de \} \in
\C_a$;

(2) the Hamiltonian flow, generated by $f$, moves $\ga$ to certain $\phi_f
(\ga)$ such that
$$
\ga \cap \phi_f(\ga) = \emptyset.
$$
For any contractible loop $\ga \in \C_a^*$ such a function can be constructed
just by hands: since $\ga$ is contractible, there exists a semi infinite curve
passing from the origin to  infinity which does not intersect $\ga$. Then this
curve together with $D^2(\de)$ can be taken as the zero level of a non negative
function. Then other level lines of such a function can go around the origin
such that the corresponding flow shall move any point of $\C_a - D^2(\de)$ from
any neighborhood of the origin to infinity, and this shows that the desired $f$
exists.

Then we can construct certain function $F$ on $\C^{k+1}$ which will generate
the flow moving $S_{\ga, (c_1, ..., c_k)}$ to $S_{\phi_f(\ga), (c_1, ...,
c_k)}$ for any couple of non critical values $(c_1, ..., c_k)$.

Let's present the explicit construction of such an $F$. Take a smaller disc
$D^2(\de_1)$ in $\C_a$ so $\de_1 < \de$, consider the subspace
$$
\psi^{-1}(\C_a - D^2(\de_1) \subset \C^{k+1},
$$
denote it as $\C^{k+1}_{\de_1}$ and the restriction $\psi$ to it as
$\psi_{\de_1}$. Then since $\psi_{\de_1}$ is smooth we have the lifted function
$\psi_{\de_1}^* f$ on $\C^{k+1}_{\de_1}$, and since this lifted function
vanishes on a neighborhood of the boundary component $\psi^{-1}(D^2(\de_1))$,
it can be continued to whole $\C^{k+1}$. We denote an appropriate continuation
by $F$, and this is our Hamiltonian.

A natural question arises: for a given function $f$ one has two vector fields
on $\C^{n+1}$: the Hamiltonian vector field of $F$ and the lifted vector field
$\nabla_{\psi} X_f$. It's not hard to see that they are proportional. It means
that the Hamiltonian flow of $F$ moves fibers to fibers outside of the singular
set. Consequently in a time this flow shall move $S_{\ga, (c_1, ..., c_k)}$ to
$S_{\phi_f \ga, (c_1, ..., c_k)}$, but since the fibers of $\psi$ do not
intersect each other it means that $F$ generates Hamiltonian isotopy which
indicates displaceability of $S_{\ga, (c_1, ..., c_k)}$.

Thus we proof the following
\begin{prop} Every torus $S_{\ga, (c_1, ..., c_k)}$ is displaceable in $\C^{k+1}$.
\end{prop}

But due to the Proposition 1.2 it implies
\begin{cor} Twist torus $\Theta^k$ is displaceable in $\C^{k+1}$.
\end{cor}

Now, the constructive proof of Proposition 2.1 hints that the same method can
be exploited in the studies on the classification problem of lagrangian tori up
to Hamiltonian isotopy. Indeed, let $\ga_1$ and $\ga_2$ be two smooth loops in
$\C_a$ which restrict domains of the same symplectic area and we do not pass
through the origin (whatever they are contractible or not --- but both of the
same type). Then there is a function $f$ which vanishes on a certain small disc
$D^2(\de) \subset \C_a$ and which moves in a finite time $\ga_1$ to $\ga_2$.
Then we can apply the same lifting construction to get a global function $F$ on
$\C^{k+1}$ whose Hamiltonian flow shall move lagrangian torus $S_{\ga_1, (c_1,
..., c_k)}$ to $S_{\ga_2, (c_1, ..., c_k)}$, so these tori are from the same
equivalence class. On the other hand, if $\ga_1$ and $\ga_2$ lie from the
different equivalence classes in $\C_a$ then $S_{\ga_1, (c_1, ..., c_k)},
S_{\ga_2, (c_1, ..., c_k)}$ must be unequivalent in $\C^{k+1}$ since they have
different Hamiltonian invariants. Thus we get

\begin{prop} Lagrangian tori $S_{\ga_1, (c_1, ..., c_k)}, S_{\ga_2, (c_1, ..., c_k)}$
in $\C^{k+1}$ both of the standard or of the Chekanov type are equivalent up to
Hamiltonian isotopy if and only if the symplectic area of the domains
restricted by $\ga_i$ is the same.
\end{prop}

A natural question arises: are there lagrangian tori of different types which
are equivalent? The construction we've used above can be modified such that one
proves the following

\begin{prop} Lagrangian tori $S_{\ga_1, (c_1, ..., c_k)}, S_{\ga_2, (c_1, ..., c_k)}$
in $\C^{k+1}$ are equivalent if

(1) $\ga_1$ and $\ga_2$ are equivalent in $\C_a$;

(2) in the couple of values $(c_1,..., c_k)$ the number of zeros is less than
2.
\end{prop}

The proof is similiar to the given one in [5]. Let's stress that this gives
nothing for the relations between the standard torus  and $\Theta^k$ in
$\C^{k+1}$ since there $c_1 = ... =  c_k = 0$.

\section{Final remarks}

At the end we would like to mention certain remarks.

First, the twist tori, constructed in [1] for $\C \PP^k$ and $S^2 \times ...
S^2$  can be described in terms of pseudo toric structures on the projective
space and on the product of projective lines. F.e., for $X = S^2_1 \times ...
\times S^2_k$ it can be done as follows. Fix homogenious coordinates $[x_j:
y_j]$ on each $S^2_j = \C \PP^1_j$ and consider the map
$$
\psi: X - B \to S^2_w = \C \PP^1_w,
$$
defined by
$$
w_0 = x_1 \cdot ... \cdot x_k, \quad w_1 = y_1 \cdot ... \cdot y_k;
$$
here $B$ is the base set of the pencil
$$
(\C \PP^1_w)^{\vee} = \{ \al x_1 \cdot ... \cdot x_k + \be y_1 \cdot ... y_k =
0, \quad  \al, \be \in \C\}
$$
in $X$. Then the compactified fiber $D_p = \overline{\psi^{-1}(p)}$ is a smooth
divisor in $X$ (except $p$ has coordinates $[0: 1]$ or $[1: 0]$ on $\C \PP^1_w$
so when $p$ is either the north or the south pole) which is a toric variety
with the same convex polytop $P$. The open part $D_p - B$ is fibered on smooth
tori, and the choice of an inner point in $P$ defines a lagrangian torus in
each $D_p$. Then for any smooth loop $\ga \subset \C \PP^1_w$ which does not
pass through the north of the south pole one has a smooth lagrangian torus in
$X$. Again we have two cases: either $\ga$ is not contractible in $\C \PP^1_w$
with punctured poles (the standard type) or $\ga$ is contractible there (the
Chekanov type).

 However for  this case Proposition 3 should not
work in full generality. Indeed, in its proof above we use the fact that a
given contractible  loop $\ga$ is displaceable in $\C_a$. Now for the compact
case it is not completely true: if a base loop $\ga$ on the base space $\C
\PP^1_w$ is {\it meridional} so if it divides $\C \PP^1_w$ into two parts of
the same symplectic area then $\ga$ is undisplaceable ifself. Therefore our
arguments work in the case of non meridional contractible loop $\ga \in \C
\PP^1_w$ only; at the same time one can mention here that this implies that our
arguments do not cover exactly the case of monotone lagrangian tori in $X$, and
I would stress that these and only these exotic tori of Chekanov and Schlenk
are Bohr --- Sommerfeld with respect to the anticanonical class.

Second, more complicated twist tori from [1] ($\Theta^k_l$ in the notation of
[1]) can be constructed in the way using slightly different pseudotoric
structures on the same space, or we can take different tori on fibers of $\psi$
instead of standard $S_{(c_1, ..., c_k)}^a$. And again the pseudo toric
considerations can help in the studies of the displaceability property.

\end{document}